\documentclass[12pt,a4paper]{article}
\usepackage{}
\usepackage{indentfirst}
\usepackage{amsfonts}
\usepackage{amssymb}
\usepackage{mathrsfs}
\usepackage{amsmath}
\usepackage{amsthm}
\usepackage{enumerate}
\usepackage{cite}
\usepackage[all]{xy}
\allowdisplaybreaks
\newtheorem{defn}{Definition}[section]
\newtheorem{thm}[defn]{Theorem}
\newtheorem{lem}[defn]{Lemma}
\newtheorem{prop}[defn]{Proposition}
\newtheorem{cor}[defn]{Corollary}
\newtheorem{eg}[defn]{Example}
\newtheorem{re}[defn]{Remark}
\newcommand{\bdefn}{\begin{defn}}
\newcommand{\edefn}{\end{defn}}
\newcommand{\bthm}{\begin{thm}}
\newcommand{\ethm}{\end{thm}}
\newcommand{\blem}{\begin{lem}}
\newcommand{\elem}{\end{lem}}
\newcommand{\bprop}{\begin{prop}}
\newcommand{\eprop}{\end{prop}}
\newcommand{\bcor}{\begin{cor}}
\newcommand{\ecor}{\end{cor}}
\newcommand{\beg}{\begin{eg}}
\newcommand{\eeg}{\end{eg}}
\newcommand{\bre}{\begin{re}}
\newcommand{\ere}{\end{re}}
\newcommand{\bpf}{\begin{proof}}
\newcommand{\epf}{\end{proof}}

\newcommand{\benu}{\begin{enumerate}}
\newcommand{\eenu}{\end{enumerate}}
\newcommand{\bc}{\begin{center}}
\newcommand{\ec}{\end{center}}
\newcommand{\bea}{\begin{eqnarray}}
\newcommand{\eea}{\end{eqnarray}}
\newcommand{\Bea}{\begin{eqnarray*}}
\newcommand{\Eea}{\end{eqnarray*}}
\newcommand{\beq}{\begin{equation}}
\newcommand{\eeq}{\end{equation}}
\newcommand{\Beq}{\begin{equation*}}
\newcommand{\Eeq}{\end{equation*}}
\newcommand{\bspl}{\begin{split}}
\newcommand{\espl}{\end{split}}
\newcommand\relphantom[1]{\mathrel{\phantom{#1}}}

\numberwithin{equation}{section}

\begin{document}
\title{\bf  Super-biderivations of the contact Lie superalgebra $K(m,n;\underline{t})$}
\author{\normalsize \bf Xiaodong Zhao$^{1,2}$, Yuan Chang$^3$, Xin Zhou$^{1,2}$, Liangyun Chen$^1$ }
\date{\small{$^{1}$ School of Mathematics and Statistics, Northeast Normal University, \\ Changchun, 130024, CHINA\\
$^{2}$ School of Mathematics and Statistics, Yili Normal University, \\ Yining, 835000, CHINA\\
$^{3}$School of Mathematics, Dongbei University of Finance and Economics,\\ Dalian, 116025, CHINA}}
\maketitle
\begin{abstract}
Let $K$ denote the contact Lie superalgebra $K(m,n;\underline{t})$ over a field of characteristic $p>3$, where $m,~n\in \mathbb{N}+1$ and $\underline{t}=(t_1,t_2,\ldots,t_m)$ is an $m$-tuple of positive integers, and $K$ has a finite $\mathbb{Z}$-graded structure.
     Let $T_K$ be the canonical torus of $K$, which is an abelian subalgebra of $K_{0}$ and semisimple on $K_{-1}$.
     Utilizing the weight space decomposition of $K$ with respect to $T_K$,
     we prove in this paper that each skew-symmetric super-biderivation of $K$ is inner.

\end{abstract}
\footnote[0]{ Corresponding author(L. Chen): chenly640@nenu.edu.cn.}
\footnote[0]{Supported by NNSF of China (Nos. 11771069), NSF of Jilin province (No. 20170101048JC) and
the project of jilin province department of education (No. JJKH20180005K).}
\section{Introduction}
Let $L$ be a Lie algebra over an arbitrary field $\mathbb{F}$.
An $\mathbb{F}$-linear map $D:L\rightarrow L$ is a derivation satisfying $$D([x,y])=[D(x),y]+[x,D(y)],$$ for all $x,y\in L$.
A bilinear map $\psi:L\times L\rightarrow L$ is called a biderivation if it is a derivation with respect to both components,
meaning that
\begin{align}
\psi(x,[y,z])&=[\psi(x,y),z]+[y,\psi(x,z)],\\
\psi([x,y],z)&=[\psi(x,z),y]+[x,\psi(y,z)],
\end{align}
for all $x,~y,~z\in L$.
A biderivation $\psi$ is called skew-symmetric if $\psi(x,y)=-\psi(y,x)$ for all $x,~y\in L$.
Obviously, if a biderivation $\psi$ is skew-symmetric, we can omit one of the equations (1.1) and (1.2).
Meanwhile, we can view $\psi(x,\cdot)$ or $\psi(\cdot,x)$ as a derivation of $L$.

The study of biderivations traces back to the research on the commuting map in the associative ring \cite{Bm}, where the author showed that all biderivations on associative prime rings are inner.
The notation of biderivations of Lie algebras was introduced in \cite{WdYxCHz}.
In recent years, there exist a lot of interests in studying biderivations and commuting maps on Lie algebras\cite{CHz,HxWdXc,Tx,WdYx,BmZHk,CHyCHly,CHyCHlyZHx}.
Moreover, the authors gave the notion of the skew-symmetric super-biderivation in \cite{XcWdHx}.
So the results about the skew-symmetric super-biderivation of Lie superalgebras arise in \cite{XcWdHx,YjTx,YchLchYc}.

The Cartan modular Lie superalgebra is an important branch of the modular Lie superalgebra, which is a Lie superalgebra over an algebraically closed field of characteristic $p>0$.
And the contact Lie superalgebra $K(m,n;\underline{t})$ is an important class of Cartan modular Lie superalgebras.
There are many research results about the contact Lie superalgebra $K(m,n;\underline{t})$, such as, derivation superalgebras\cite{FmQzh,BgWl,BgLch}, noncontractible filtrations\cite{YZhWl}, nondegenerate associative bilinear forms\cite{YwYZh}.

%In this paper, we use the weight space decomposition of $K(m,n;\underline{t})$ with respect to the canonical torus $T_K$ and prove all skew-symmetric super-biderivations of the contact Lie superalgebra $K(m,n;\underline{t})$ are inner.
%%
%The paper is organized as follows.
%%
%In Section 2, we review the basic notation.
%%
%In Section 3, utilizing the canonical torus $T_K$ and the weight space decomposition of $K(m,n;\underline{t})$ with respect to $T_K$,
%we show each derivation composed of the skew-symmetric super-biderivation and the element of $T_K$ is inner.
%Moreover, we prove that each skew-symmetric super-biderivation of $K(m,n;\underline{t})$ is inner.
In this paper, we prove that each skew-symmetric super-biderivation of $K(m,n;\underline{t})$ is inner.
The paper is organized as follows.
In Section 2, we recall the basic notation.
In Section 3, we use the weight space decomposition of $K(m,n;\underline{t})$ with respect to the canonical torus $T_K$ to prove that all skew-symmetric super-biderivation of $K$ is inner(Theorem 3.14).

\section{Preliminaries}

Let $\mathbb{F}$ denote the prime field of the characteristic $p>2$ and $\mathbb{Z}_{2}=\{\overline{0},\overline{1}\}$ the additive group of two elements.
For a vector superspace $V=V_{\overline{0}}\oplus V_{\overline{1}}$, we use $\mathrm{p}(x)$ for the parity of $x\in V_{\alpha}$, $\alpha\in\mathbb{Z}_{2}$.
If $V=\oplus_{i\in\mathbb{Z}}V_i$ is a $\mathbb{Z}$-graded vector space and $x\in V$ is a $\mathbb{Z}$-homogeneous element, write $|x|$ for the $\mathbb{Z}$-degree of $x$.
Once the symbol $\mathrm{p}(x)$ or $|x|$ appears in this paper, it implies that $x$ is a $\mathbb{Z}_{2}$-homogeneous element or that $x$ is a $\mathbb{Z}$-homogeneous element.
Throughout this paper all vector spaces or algebras are over $\mathbb{F}$.
\subsection{Skew-symmetric super-biderivations of a Lie superalgebra}
Let us recall some facts related to the superderivation and skew-symmetric super-biderivation of Lie superalgebras.
A Lie superalgebra is a vector superspace $L=L_{\overline{0}}\oplus L_{\overline{1}}$ with an even bilinear mapping  $[\cdot,\cdot]:L\times L\rightarrow L$ satisfying the following axioms:
\begin{align*}
[x,y]&=-(-1)^{\mathrm{p}(x)\mathrm{p}(y)}[y,x],\\
[x,[y,z]]&=[[x,y],z]+(-1)^{\mathrm{p}(x)\mathrm{p}(y)}[y,[x,z]],
\end{align*}
for all $x,y,z\in L$.
We call a linear mapping $D:L\times L\rightarrow L$ a superderivation of $L$ if it satisfies the following axiom:
$$D([x,y])=[D(x),y]+(-1)^{\mathrm{p}(D)\mathrm{p}(x)}[x,D(y)],$$
for all $x,y\in L$, where $\mathrm{p}(D)$ denotes the $\mathbb{Z}_{2}$-degree of $D$.
Write $\mathrm{Der}_{\overline{0}}(L)$ (resp. $\mathrm{Der}_{\overline{1}}(L)$) for the set of all superderivations of $\mathbb{Z}_2$-degree $\overline{0}$ (resp. $\overline{1}$) of $L$.

We call a bilinear mapping  $\phi:L\times L\rightarrow L$
a skew-symmetric super-biderivation of $L$ if it satisfies the following axioms:
\begin{align*}
&skew-symmetry:~~\phi(x,y)=-(-1)^{\mathrm{p}(x)\mathrm{p}(y)}\phi(y,x),\\
&\phi([x,y],z)=(-1)^{\mathrm{p}(\phi)\mathrm{p}(x)}[x,\phi(y,z)]+(-1)^{\mathrm{p}(y)\mathrm{p}(z)}[\phi(x,z),y],\\
&\phi(x,[y,z])=[\phi(x,y),z]+(-1)^{(\mathrm{p}(\phi)+\mathrm{p}(x))\mathrm{p}(y)}[y,\phi(x,z)],
\end{align*}
for all $\mathbb{Z}_2$-homogeneous elements $x,y,z\in L$.
A super-biderivation $\phi$ of $\mathbb{Z}_{2}$-degree $\gamma$  of $L$ is a super-biderivation such that $\phi(L_{\alpha},L_{\beta})\subset L_{\alpha+\beta+\gamma}$ for any $\alpha,~\beta \in \mathbb{Z}_{2} $.
Denote by $\mathrm{BDer}_{\gamma}(L)$ the set of all skew-symmetric super-biderivations of $\mathbb{Z}_2$-degree $\gamma$. Obviously,
$$\mathrm{BDer}(L)=\mathrm{BDer}_{\overline{0}}(L)\oplus \mathrm{BDer}_{\overline{1}}(L).$$
Specially, if the bilinear map $\phi_{\lambda}:L\times L\rightarrow L$ is defined by $\phi_{\lambda}(x,y)=\lambda[x,y]$ for all $x,y\in L$, where $\lambda\in \mathbb{F}$, then it is easy to check that $\phi_{\lambda}$ is a super-biderivation of $L$.
This class of super-biderivations is called inner.
Denote by  $\mathrm{IBDer}(L)$ the set of all inner super-biderivations.
\subsection{Contact Lie superalgebras $K(m,n;\underline{t})$}
We propose to construct a $\mathbb{Z}_2$-gradation tensor algebra via a divided~power~algebra and a exterior~superalgebra.
In the follow, we introduce the divided~power~algebra $\mathcal{O}(m)$ and the exterior~superalgebra $\Lambda(n)$.
Fix two positive integers $m>1$ and $n>1$.
For $\alpha=(\alpha_1,\ldots,\alpha_m)\in\mathbb{N}^{m}$, where $ \mathbb{N}$ denote the set of natural numbers, put $|\alpha|=\Sigma_{i=1}^{m}\alpha_i$.
For two $m$-tuples $\alpha=(\alpha_1,\ldots,\alpha_m)$ and $\beta=(\beta_1,\ldots,\beta_m)\in\mathbb{N}^m$,
we write ${\alpha\choose \beta}=\prod_{i=1}^{m}{\alpha_i\choose \beta_i}$
and define $\beta\leq \alpha \Longleftrightarrow \beta_i\leq \alpha_i$, $1\leq i\leq m$.
Let $\mathcal{O}(m)$ denote the $\mathbb{F}$-algebra of divided power series in the variable $x_1,\ldots,x_m$, which is called a $divided~power~algebra$.
For convenience, we replace $x_1^{\alpha_1}x_2^{\alpha_2}\cdots x_m^{\alpha_m}$ by $x^{(\alpha)}$, $\alpha=(\alpha_1,\alpha_2,\ldots,\alpha_m)$.
Obviously, $\mathcal{O}(m)$ has an $\mathbb{F}$-basis $\{x^{(\alpha)}|\alpha\in\mathbb{N}^{m}\}$ and satisfies the formula:
\begin{eqnarray}\label{O¼ÆË㹫ʽ1}
x^{(\alpha)}x^{(\beta)}={\alpha+\beta\choose\alpha}x^{(\alpha+\beta)},~\forall~\alpha,\beta\in\mathbb{N}^{m}.
\end{eqnarray}
Let $\Lambda(n)$ denote the $exterior~superalgebra$ over $\mathbb{F}$ with $n$ variables $x_{m+1},\ldots,x_{s}$, where $s=m+n$.
The tensor product $\mathcal{O}(m,n)=\mathcal{O}(m)\otimes_{\mathbb{F}}\Lambda(n)$ is an associative superalgebra with a $\mathbb{Z}_2$-gradation induced by the trivial $\mathbb{Z}_2$-gradation of $\mathcal{O}(m)$ and the natural $\mathbb{Z}_2$-gradation of $\Lambda(n)$.
Obviously,~$\mathcal{O}(m,n)$ is super-commutative.
For $g\in\mathcal{O}(m)$,~$f\in\Lambda(n)$,~it is customary to write $gf$ instead of $g\otimes f$.
Including the formula (\ref{O¼ÆË㹫ʽ1}), the following formulas also hold in $\mathcal{O}(m,n)$:
$$x_kx_l=-x_lx_k,~\forall~k,l\in\{m+1,\ldots,s\};$$
$$x^{(\alpha)}x_k=x_kx^{(\alpha)},~\forall~\alpha\in\mathbb{N}^{m}, k\in\{m+1,\ldots,s\}.$$
For $k=1,\ldots,n$, set
$$\mathbb{B}_k:=\{\langle i_1,i_2,\ldots,i_k\rangle\mid m+1\leq i_1<i_2<\ldots<i_k\leq s\}$$
and~$\mathbb{B}:=\cup_{k=0}^{n}\mathbb{B}_k$,~where~$\mathbb{B}_0=\emptyset$.
For~$u=\langle i_1,i_2,\ldots,i_k\rangle\in\mathbb{B}_k$,~set~$|u|:=k$,~$x^u=x_{i_1}\cdots x_{i_k}$.
Specially, we define~$|\emptyset|=0$, $x^{\emptyset}=1$, $|\omega|=n$ and $x^{\omega}=x_{m+1}\cdots x_{m+n}$.
Clearly,~the set $\{x^{(\alpha)}x^u|~\alpha\in\mathbb{N}^{m},u\in\mathbb{B}\}$~constitutes an $\mathbb{F}$-basis of $\mathcal{O}(m,n)$.

Put~$\mathrm{I}_0:=\{1,\ldots,m\}$,~$\mathrm{I}_1:=\{m+1,\ldots,m+n\}$~and~$\mathrm{I}:=\mathrm{I}_0\cup \mathrm{I}_1$.
Let~$\partial_1,\partial_2,\ldots,\partial_{s}$ be the linear transformations of~$\mathcal{O}(m,n)$ such that $\partial_i(x^{(\alpha)})=x^{(\alpha-\varepsilon_i)}$ for $i\in \mathrm{I}_0$, %$\alpha\in\mathbb{N}^{m}$
and $\partial_i(x_k)=\delta_{ik}$, $k\in \mathrm{I}_1$, for $i\in \mathrm{I}_1$, where $\delta_{ij}$ is denoted the Kronecker symbol.
Obviously,~$\mathrm{p}(\partial_i)=\overline{0}$ if $i\in \mathrm{I}_0$ and~$\mathrm{p}(\partial_i)=\overline{1}$ if $i\in \mathrm{I}_1$.
Then $\partial_1,\partial_2,\ldots,\partial_{s}$ are superderivations of the superalgebra $\mathcal{O}(m,n)$.
Let
$$W(m,n):=\left\{\sum f_r\partial_r|~f_r\in\mathcal{O}(m,n),r\in\mathrm{I}\right\}.$$%_{r\in \mathrm{I}}
Then $W(m,n)$ is an infinite-dimensional Lie superalgebra contained in $\mathrm{Der}(\mathcal{O}(m,n))$.
One can verify that
\begin{eqnarray}\label{W()¼ÆË㹫ʽ}
[f\partial_i,g\partial_j]=f\partial_i(g)\partial_j-(-1)^{\mathrm{p}(f\partial_i)\mathrm{p}(g\partial_j)}g\partial_j(f)\partial_i,
\end{eqnarray}
for all $f,g\in\mathcal{O}(m,n)$ and $i,j\in \mathrm{I}$.

Fix two $m$-tuples of positive integers $\underline{t}=(t_1,t_2,\ldots,t_m)$ and $\pi=(\pi_1,\pi_2,\ldots,\pi_m)$, where $\pi_i=p^{t_i}-1$ for all $i\in\mathrm{I}_0$ and $p$ is denoted the characteristic of the basic field $\mathbb{F}$.
For two $m$-tuples $\alpha=(\alpha_1,\ldots,\alpha_m)$ and $\beta=(\beta_1,\ldots,\beta_m)\in\mathbb{N}^m$,
we have ${\alpha+\beta\choose \alpha}=0$ if there is some $i\in\{1,\ldots,m\}$ satisfying $\alpha_i+\beta_i\geq p^{t_i}$.
Thence the set $$O(m,n;\underline{t})=\{x^{(\alpha)}x^u~|~0\leq \alpha \leq \pi, u\in\mathbb{B}\}$$ is a subalgebra of $O(m,n)$
and the set $$W(m,n;\underline{t})=\mathrm{span}_{\mathbb{F}}\{x^{(\alpha)}x^u\partial_r~|~0\leq \alpha \leq \pi, u\in\mathbb{B}, r\in \mathrm{I}\}$$
is a finite-dimensional simple subalgebra of $W(m,n)$, which is called the generalized Witt Lie superalgebra.
$W(m,n;\underline{t})$ possesses a $\mathbb{Z}$-graded structure:
$$W(m,n;\underline{t})=\bigoplus_{r=-1}^{\xi-1}W(m,n;\underline{t})_r,$$
where $W(m,n;\underline{t})_r:=\mathrm{span}_{\mathbb{F}}\{x^{(\alpha)}x^u\partial_j|~|\alpha|+|u|=r+1,~j\in\mathrm{I}\}$ and $\xi:=|\pi|+n$.
For $i\in \mathrm{I}_0$, we abbreviate $x^{(\varepsilon_i)}$ to $x_i$, where $\varepsilon_i$ is denoted the $m$-tuple with 1 as the i-th entry and 0 elsewhere.

Hereafter, suppose $m=2r+1$ is odd and $n=2t$ is even.
Let $\mathrm{J}=\mathrm{I} \setminus\{m\}$ and $\mathrm{J}_0=\mathrm{I}_0 \setminus\{m\}$.
For $i\in \mathrm{J}$, put\\\[i':=\begin{cases}i+r,&1\leq i\leq r,\\i-r,  &r<i\leq 2r,\\i, &i=m,\\i+t, &m<i\leq m+t,\\i-t, &m+t<i\leq s;
\end{cases}~~~~\sigma(i):=\begin{cases}1, &1\leq i\leq r\\-1,  &r<i\leq 2r\\1, &2r<i\leq s.
\end{cases}\]
Define a linear mapping $D_{K}:\mathcal{O}(m,n)\rightarrow W(m,n)$ by means of
$$D_{K}(f)=\sum \limits_{i\in \mathrm{J}}(-1)^{\mathrm{p}(\partial_i)\mathrm{p}(f)}(x_{i}\partial_{m}(f)+\sigma(i')\partial_{i'}(f))\partial_{i}+
(2f-\sum \limits_{i\in \mathrm{J}}x_{i}\partial_{m}(f))\partial_{m}.$$
%
%Then $K(m,n)=\left\{D_{k}(f)|~f\in\mathcal{O}(m,n)\right\}$ is an infinite-dimensional subsuperalgebra of W(m,n). Put
%%
%$$K(m,n;\underline{t})=\left\{D_{k}(f)|~f\in\mathcal{O}(m,n;\underline{t})\right\}$$
%
%Here ${O}(m,n;\underline{t})$ is a finite dimensional subsuperalgebra of ${O}(m,n)$ with a  $\mathbb{Z}$-graded structure:
%$${O}(m,n;\underline{t})=\bigoplus_{r=-1}^{\xi-1}{O}(m,n;\underline{t})_r,$$
%
The restricted linear mapping of $D_{K}$ on $\mathcal O(m,n;\underline{t})$ still is denoted by $D_{K}$, that is
$$D_{K}:\mathcal{O}(m,n;\underline{t})\rightarrow K(m,n;\underline{t}).$$
Let $\widetilde{K}(m,n;\underline{t})$ denote the image of $\mathcal O(m,n;\underline{t})$ under $D_{K}$.
Consider the derived algebra of $\widetilde{K}(m,n;\underline{t})$:$$K(m,n;\underline{t})=[\widetilde{K}(m,n;\underline{t}),\widetilde{K}(m,n;\underline{t})].$$
The derived algebra $K(m,n;\underline{t})$ is a finite dimensional simple Lie superalgebra, which is called the contact Lie superalgebra.
We define a Lie bracket $\langle\cdot,\cdot\rangle$ on the tensor superalgebra $\mathcal{O}(m,n;\underline{t})$ by
$$ \langle f, g \rangle:=D_{K}(f)(g)-2\partial_{m}(f)(g),$$
for all $f,g\in\mathcal{O}(m,n;\underline{t})$.
Since $D_{K}$ is injective and $D_{K}(\langle f,g\rangle)=[ D_{K}(f),D_{K}(g)]$,
there exists an isomorphism, that is, $$({K}(m,n;\underline{t}),[\cdot,\cdot])\cong(\mathcal{O}(m,n;\underline{t}),\langle\cdot,\cdot\rangle).$$
For convenience, we use $W$ and $W_r$ denote $W(m,n;\underline{t})$ and its $\mathbb{Z}$-graded subspace $W(m,n;\underline{t})_r$, respectively, ${K}(m,n;\underline{t})$ is denoted by $K$.

\section{Skew-symmetry Super-biderivation of $K(m,n;\underline{t})$}

\begin{lem}\cite{FgDx}\label{б³¬Ë«1}
Let $L$ be a Lie superalgebra.
Suppose that $\phi$ is a skew-symmetric super-biderivation on $L$, then $$[\phi(x,y),[u,v]]=(-1)^{\mathrm{p}(\phi)(\mathrm{p}(x)+\mathrm{p}(y))}[[x,y] ,\phi(u,v)] $$ for any homogenous element $x,y,u,v\in L$.
\end{lem}

\begin{lem}\cite{FgDx}\label{б³¬Ë«2}
Let $L$ be a Lie superalgebra.
Suppose that $\phi$ is a skew-symmetric super-biderivation on $L$.
If $\mathrm{p}(x)+\mathrm{p}(y)=\overline{0}$, then $$[\phi(x,y),[x,y]]=0$$ for any homogenous element $x,y\in L$.
\end{lem}

\begin{lem}\cite{FgDx}\label{б³¬Ë«3}
Let $L$ be a Lie superalgebra.
Suppose that $\phi$ is a skew-symmetric super-biderivation on $L$. If $[x,y]=0$, then $\phi(x,y)\in C_L([L,L])$, where $C_L([L,L])$ is the centralizer of $[L,L]$ .
\end{lem}
\begin{lem}\label{Kб³¬Ë«1}
Let $K$ denote the contact Lie superalgebra.
Suppose $\phi$ is a skew-symmetric super-biderivation on $K$. If $\langle x,y \rangle=0$ for $x, y\in K$,  then $\phi(x,y)=0 .$
\end{lem}
\bpf
Since $K$ is a simple Lie superalgebra, it is obvious that $K=\langle K,K \rangle $ and $C(K)=0$.
if $\langle x,y \rangle=0$ for $x, y\in K$, by Lemma \ref{б³¬Ë«3}, we obtain $\phi(x,y)\in C_K(\langle K,K \rangle)=C(K)=0.$
\epf

Set $T_K=\mathrm{span}_{\mathbb{F}}\{x_{i}x_{i'}~|~i\in\mathrm{J}\}$.
Obviously, $T_{K}\subseteq K(m,n;\underline{t})_{0}\bigcap K(m,n;\underline{t})_{\bar{0}}$.
$T_K$ is an abelian subalgebra of $K$.
For any $x^{(\alpha)}x^u\in K$, we have
\begin{eqnarray}
\langle x_{i}x_{i'},x^{(\alpha)}x^u \rangle=(\alpha_{i'}-\alpha_{i}+\delta_{(i'\in u)}-\delta_{(i\in u)})x^{(\alpha)}x^u,
\end{eqnarray}
where $\delta_{(\mathrm{P})}=1$ if the proposition $\mathrm{P}$ is true, $=0$  if the proposition $\mathrm{P}$ is false.
Fixed an $m$-tuple $\alpha$, where  $\alpha\in \mathbb{N}^m$, $0\leq\alpha \leq \pi$ and $u\in \mathbb{B}$,~we define a linear function $(\alpha+\langle u\rangle): T\rightarrow \mathbb{F}$ such that
$$(\alpha+\langle u\rangle)(x_{i}x_{i'})= \alpha_{i'}-\alpha_{i}+\delta_{(i'\in u)}-\delta_{(i\in u)}.$$
Further, $K$ has a weight space decomposition with respect to $T_K$:
$$K=\bigoplus_{(\alpha+\langle u\rangle)}K_{(\alpha+\langle u\rangle)}.$$
%where
% $$K_{(\alpha+\langle u\rangle)}=\mathrm{span}_{\mathbb{F}}\{x^{(\alpha)}x^u\in O| \langle x_{i}x_{i'},x^{(\alpha)}x^u \rangle=(\alpha_{i'}-\alpha_{i}+\delta_{(i'\in u)}-\delta_{(i\in u)})x^{(\alpha)}x^u, \\\forall x_{i}x_{i'}\in T\} .$$
\begin{lem}\label{TµÄб˫³¬1}
Suppose that $\phi$ is a $\mathbb{Z}_2$-homogeneous skew-symmetric super-biderivation on $K$.
Let $x^{(\alpha)}x^u\in K$ such that
$$\phi(x_{i}x_{i'},x^{(\alpha)}x^u)\in K_{(\alpha+\langle u\rangle)},$$
for any $x_{i}x_{i'}\in T_K$.
\end{lem}
\bpf
The equation by Lemma \ref{Kб³¬Ë«1},
it follows that $\phi(x_{i}x_{i'},x_{j}x_{j'})=0$ for any $i,~j\in \mathrm{J}$ from  $[x_{i}x_{i'},x_{j}x_{j'}]=0$.
Note that $\mathrm{p}(x_{l}x_{l'})=\overline{0}$ for all $l\in \mathrm{J}$, then all $x^{(\alpha)}x^u\in K$, it is clear that
\begin{align*}
&\langle x_{l}x_{l'},\phi(x_{i}x_{i'},x^{(\alpha)}x^u) \rangle\\
=&(-1)^{(\mathrm{p}(\phi)+\mathrm{p}(x_{i}x_{i'}))\mathrm{p}(x_{l}x_{l'})}(\phi(x_{i}x_{i'},\langle x_{l}x_{l'},x^{(\alpha)}x^u \rangle)-\langle \phi(x_{i}x_{i'},x_{l}x_{l'}),x^{(\alpha)}x^u \rangle)\\
=&(\alpha_{i'}-\alpha_{i}+\delta_{(i'\in u)}-\delta_{(i\in u)})\phi(x_{i}x_{i'},x^{(\alpha)}x^u).
\end{align*}
The proof is completed.
\epf

\begin{re}\label{б³¬Ë«ÊÇżµÄ}
Due to Lemma \ref{TµÄб˫³¬1}, we can find that any $\mathbb{Z}_2$-homogeneous skew-symmetric super-biderivation on $K$ is an even bilinear map.
Since $\phi(x_{i}x_{i'},x^{(\alpha)}x^u)$ and $x^{(\alpha)}x^u$ have the same $\mathbb{Z}_2$-degree.
Then the $\mathbb{Z}_2$-degree of $\phi$ is even.
\end{re}

\begin{lem}\cite{YZhWl}\label{kµÄÉú³ÉÔª}
Let $M=\{x^{(\kappa_i\varepsilon_i)}|0\leq \kappa_i\leq\pi_i, i\in \mathrm{I}_0\}$ and $N=\{x_i|i\in \mathrm{I}_1\}$. Then $K$ is generated by $M\cup N$.
\end{lem}

%
%We want to describe the action of  the skew-symmetric super-biderivation $\phi$ on $K$ if only and if we can give the action of $\phi$ on the generators of $K$.
%We give the weight spaces with the same weight of the generators of $K$ about the weight space decompositions with respect to $T$.

\blem\label{WÉú³ÉÔªµÄȨ¿Õ¼ä}
Let $i\in \mathrm{J}_0$, $j \in \mathrm{I}_1$
and $q_i\in\mathbb{N}$, $1\leq q_i\leq \pi_i$. Then the following statements hold:
\begin{itemize}
%\item[\rm{(1)}]  $K(m,n;\underline{t})_{(\varepsilon_{i})}
%=\sum\limits_{0\leq\alpha\leq\pi}\mathbb{F}(\mathop{\prod\limits_{l\in\mathrm{I}_0\backslash\{i,i^{'},m\}}} \limits_{\alpha_{l^{'}}-\alpha_{l}\equiv0(modp)} x^{(\alpha_{l}\varepsilon_l)})x^{((\alpha_i-1)\varepsilon_{i'})}x^{(\alpha_i\varepsilon_i)}x^{(\alpha_m\varepsilon_m)}x^{\bar{u}};$
\item[\rm{(1)}] $K(m,n;\underline{t})_{(0)}=\sum \limits_{0\leq\alpha\leq\pi,~\bar{u}\in\mathbb{B}}\mathbb{F}(\mathop{\prod \limits_{l\in \mathrm{J}_{0}}} \limits_{\alpha_{l'}-\alpha_{l}\equiv0~(\mathrm{mod}~p)} x^{(\alpha_{l}\varepsilon_{l})})x^{(\alpha_m\varepsilon_m)}x^{\bar{u}};$
\item[\rm{(2)}]$K(m,n;\underline{t})_{(q_i\varepsilon_i)}=\sum\limits_{0\leq\alpha\leq\pi,~\bar{u}\in\mathbb{B}}\mathbb{F}(\mathop{\prod\limits_{l\in\mathrm{J}_0\backslash\{i,{i'}\}}} \limits_{\alpha_{{l'}}-\alpha_{l}\equiv0~(\mathrm{mod}~p)} x^{(\alpha_{l}\varepsilon_l)})x^{(\alpha_i\varepsilon_i)}x^{(\overline{(\alpha_i-q_{i})}\varepsilon_{i'})}x^{(\alpha_m\varepsilon_m)}x^{\bar{u}};$
%\item[\rm{(3)}] $K(m,n;\underline{t})_{(m)}=\sum \limits_{0\leq\alpha\leq\pi}\mathbb{F}(\mathop{\prod \limits_{i\in I_{0}}} \limits_{\alpha_{i^{'}}-\alpha_{i}\equiv0(modp)} x^{(\alpha_{i}\varepsilon_{i})})x^{\bar{u}};$
\item[\rm{(3)}] $K(m,n;\underline{t})_{(\langle j\rangle)}=\sum \limits_{0\leq\alpha\leq\pi,~\bar{u}\in\mathbb{B}}\mathbb{F}(\mathop{\prod \limits_{l\in \mathrm{J}_{0}}} \limits_{\alpha_{{l'}}-\alpha_{l}\equiv0~(\mathrm{mod}~p)} x^{(\alpha_{l}\varepsilon_{l})})x^{(\alpha_m\varepsilon_m)}x_{j}x^{\bar{u}};$

\end{itemize}
where $i$ and ${i'}$ are both in $\bar{u}$ for $i\in \mathrm{J}$, and $\alpha_{l}^{\overline{q}}$ is denoted some integer and $\alpha_{l}^{\overline{q}}\equiv q~(\mathrm{mod}~p)$.
\elem
\bpf
(1)
We first discuss the vector of the same weight with 1 in $K$ with respect to $T_K$. Since we have the equation
$$\left \langle{x_{l}x_{l'}},{1} \right \rangle=D_{K}(x_{l}x_{l'})(1)-2\partial_{m}(x_{l}x_{l'})(1)=0.$$
For any $l \in \mathrm{J}$, in contrast with equation~(3.1), we get that
$$\alpha_{l'}-\alpha_{l}+\delta_{(l'\in u)}-\delta_{(l\in u)}=0.$$
 Then if $l \in \mathrm{J_{0}}, $ it is obvious that is $\alpha_{l'}-\alpha_{l}\equiv0~(\mathrm{mod}~p)$.
 If $l\in \mathrm{I_{1}}$, it is obvious that $l$~and~$l^{'}$ are both in $\bar{u}$.
 It proves that
 $$K(m,n;\underline{t})_{(0)}=\sum \limits_{0\leq\alpha\leq\pi,~\bar{u}\in\mathbb{B}}\mathbb{F}(\mathop{\prod \limits_{l\in \mathrm{J_{0}}}} \limits_{\alpha_{{l'}}-\alpha_{l}\equiv0(\mathrm{mod}~p)} x^{(\alpha_{l}\varepsilon_{l})})x^{(\alpha_m\varepsilon_m)}x^{\bar{u}}.$$

(2)
Without loss of generality,~we choose a fixed element $i\in \mathrm{J}$.
For any $l \in \mathrm{J}$, we have the equation
$$\left \langle{x_{l}x_{l'}},{x^{(q_{i}\varepsilon_{i})}} \right \rangle=D_{K}(x_{l}x_{l'})(x_{i})-2\partial_{m}(x_{l}x_{l'})(x_{i})=-q_{l}x_{l}\delta_{(li)}.$$
For any $l \in \mathrm{J}$, by equation(3.1) we have that
$$\alpha_{l'}-\alpha_{l}+\delta_{(l'\in u)}-\delta_{(l\in u)}=-q_{l}\delta_{(li)}.$$
Then we try to discuss the choice of $l\in \mathrm{J}.$
If $l\in \mathrm{J_{0}}\setminus\{i,i^{'}\}$, it is obvious that $\alpha_{l'}-\alpha_{l}\equiv0~(\mathrm{mod}~p)$.
If $ l=i$, it is obvious that $\alpha_{l'}-\alpha_{l}\equiv-q_{i}~(\mathrm{mod}~p)$.
If $l \in \mathrm{I_{1}}$, we have that $l$~and~$l'$ are both in $\bar{u}$.
So we proves that
 $$K(m,n;\underline{t})_{(q_i\varepsilon_i)}=\sum\limits_{0\leq\alpha\leq\pi,~\bar{u}\in\mathbb{B}}\mathbb{F}(\mathop{\prod\limits_{l\in\mathrm{J}_0\backslash\{i,{i'}\}}} \limits_{\alpha_{{l'}}-\alpha_{l}\equiv0~(\mathrm{mod}~p)} x^{(\alpha_{l}\varepsilon_l)})x^{(\alpha_i\varepsilon_i)}x^{(\overline{(\alpha_i-q_{i})}\varepsilon_{i'})}x^{(\alpha_m\varepsilon_m)}x^{\bar{u}}.$$

(3)
Without loss of generality,~we choose a fixed element $i\in \mathrm{J}$.
For any $l \in \mathrm{J}$, we have the equation
$$\left \langle{x_{l}x_{l'}},{x_{j}} \right \rangle=D_{K}(x_{l}x_{l'})(x_{j})-2\partial_{m}(x_{l}x_{l'})(x_{j})=-x_{l}\delta_{(lj)}.$$
By equation (3.1), for any $l \in \mathrm{J}$, we have that
$$\alpha_{l'}-\alpha_{l}+\delta_{(l'\in u)}-\delta_{(l\in u)}=-\delta_{(lj)}.$$
Then we try to discuss the choice of $l\in \mathrm{J}.$
If $l \in \mathrm{J_{0}}$, it is obvious that $\alpha_{l'}-\alpha_{l}\equiv0~(\mathrm{mod}~p)$.
If $l \in \mathrm{I_{1}}$, we have that $l$~and~$l'$ are both in $\bar{u}$.
It proves that
 $$K(m,n;\underline{t})_{(\langle j\rangle)}=\sum \limits_{0\leq\alpha\leq\pi,~\bar{u}\in\mathbb{B}}\mathbb{F}(\mathop{\prod \limits_{l\in \mathrm{J_{0}}}} \limits_{\alpha_{{l'}}-\alpha_{l}\equiv0~(\mathrm{mod}~p)} x^{(\alpha_{l}\varepsilon_{l})})x^{(\alpha_m\varepsilon_m)}x_{j}x^{\bar{u}}.$$
\epf

\blem\label{3.09}
Suppose that $\phi$ is a  $\mathbb{Z}_2$-homogeneous skew-symmetric super-biderivation on $K$.
 For any $x_{l}x_{l'}\in T_K$ and $x^{(q_{m}\varepsilon_{m})}\in M$, where $0\leq q_m\leq \pi_m$, we have
%there is an element $ \lambda\in\mathbb{F}$ such that
$$\phi(x_{l}x_{l'},x^{(q_{m}\varepsilon_{m})})=0.$$
%\lambda\langle x_{l}x_{l'},x^{(q_{m}\varepsilon_{m})} \rangle
\elem

\bpf
When $q_{m}=0$, by Lemma \ref{Kб³¬Ë«1}, it is obvious that $\phi(x_{l}x_{l'},1)=0$ for $l\in\mathrm{I} $ from the equation $\langle x_{l}x_{l'},1\rangle=0$.
When $q_{m}\neq0$, we have that
 \begin{align*}
&\langle x_{l}x_{l'},x^{(q_m\varepsilon_m)} \rangle\\
=&D_{K}(x_{l}x_{l'})(x^{(q_m\varepsilon_m)})-2\partial_{m}(x_{l}x_{l'})(x^{(q_m\varepsilon_m)})\\
%=&(x_{l}x_{l'})_{m}D_{m}(x^{(q_m\varepsilon_m)})\\
=&(2(x_{l}x_{l'})-\sum\limits_{i \in \mathrm{J}}x_{i}\partial_{i}(x_{l}x_{l'}))\partial_{m}(x^{(q_m\varepsilon_m)})\\
=&(2(x_{l}x_{l'})-2(x_{l}x_{l'}))\partial_{m}(x^{(q_m\varepsilon_m)})\\
=&0.
\end{align*}
Hence, we have that
$$\phi(x_{l}x_{l'},x^{(q_{m}\varepsilon_{m})})=0.$$
%where $\lambda\in\mathbb{F}$.
The proof is completed.
\epf

\blem\label{3.10}
Suppose that $\phi$ is a  $\mathbb{Z}_2$-homogeneous skew-symmetric super-biderivation on $K$.
For any $i\in \mathrm{J} $ and $x_{l}x_{l'}\in T_K$, there is an element $ \lambda_i\in\mathbb{F}$ such that
$$\phi(x_{l}x_{l'},x_{i})=\lambda_i\langle x_{l}x_{l'},x_{i} \rangle,$$
where $\lambda_i$ is dependent on the second component.
\elem

\bpf
Without loss of generality,~we choose a fixed element $i\in \mathrm{J}$.
By Lemma \ref{Kб³¬Ë«1}, it is obvious that  $\phi(x_{l}x_{l'},x_{i})=0$ for $l\in\mathrm{J}\setminus\{i,{i'}\} $ from $\langle x_{l}x_{l'},x_{i}\rangle=0$.
So we only need to discuss the case with the condition $l=i$.

When $q_{i}=1$, by Lemma \ref{Kб³¬Ë«1}~(2), we can suppose that
\begin{align*}
\phi(x_{l}x_{l'},x_{i})=&\sum\limits_{0\leq\alpha\leq\pi,~\overline{u}\in\mathbb{B}}a(\alpha,\bar{u})(\mathop{\prod\limits_{l\in\mathrm{J}_0\backslash\{i,{i'}\}}} \limits_{\alpha_{l'}-\alpha_{l}\equiv0~(\mathrm{mod}~p)} x^{(\alpha_{l}\varepsilon_l)})x^{(\alpha_i\varepsilon_i)}x^{((\overline{\alpha_i-1})\varepsilon_{i'})}x^{(\alpha_m\varepsilon_m)}x^{\bar{u}}.
\end{align*}
It is obvious that% $\phi(x_{i}x_{i'},1)=0$, for $\left \langle{x_{i}x_{i'}},{1} \right \rangle=0$\\
\begin{align*}
0=&\phi(x_{i}x_{i'},\langle 1,x_{i} \rangle)-\langle \phi(x_{i}x_{i'},1),x_{i}\rangle\\
=&(-1)^{(\mathrm{p}(\phi)+\mathrm{p}(x_{i}x_{i'}))\mathrm{p}(1)}\langle 1,\phi(x_{i}x_{i'},x_{i}) \rangle\\
%=&\langle 1,\phi(x_{i}x_{i'},x_{i}) \rangle\\
=&D_{K}(1)( \phi(x_{i}x_{i'},x_{i}))-2\partial_{m}(1)( \phi(x_{i}x_{i'},x_{i}))\\
=&\partial_{m}( \phi(x_{i}x_{i'},x_{i}))\\
=&\partial_{m}(\sum\limits_{0\leq\alpha\leq\pi,~\overline{u}\in\mathbb{B}}a(\alpha,\bar{u})(\mathop{\prod\limits_{l\in\mathrm{J}_0\backslash\{i,{i'}\}}} \limits_{\alpha_{{l'}}-\alpha_{l}\equiv0~(\mathrm{mod}~p)} x^{(\alpha_{l}\varepsilon_l)})x^{(\alpha_i\varepsilon_i)}x^{((\overline{\alpha_i-1})\varepsilon_{i'})}x^{(\alpha_m\varepsilon_m)}x^{\bar{u}}).
\end{align*}
By computing the equation, we find that $a(\alpha,\bar{u})=0$ if $\alpha_{m}>0$.
Putting $l\in\mathrm{J}\backslash\{i,i^{\prime}\}$, we have that
%$k=l$
\begin{align*}
0=&(-1)^{(\mathrm{p}(\phi)+\mathrm{p}(x_{i}x_{i'}))\mathrm{p}(x_{l})}(\phi(x_{i}x_{i'}, \langle x_{l},x_{i} \rangle)-\langle \phi(x_{i}x_{i'},x_{l}),x_{i} \rangle)\\
=&\langle x_{l},\phi(x_{i}x_{i'},x_{i}) \rangle\\
=&D_{K}(x_{l})( \phi(x_{i}x_{i'},x_{i}))-2\partial_{m}(x_{l})( \phi(x_{i}x_{i'},x_{i}))\\
=&\partial_{l'}( \phi(x_{i}x_{i'},x_{i}))\\
=&\partial_{l'}(\sum\limits_{0\leq\alpha\leq\pi,~\overline{u}\in\mathbb{B}} a(\alpha,i)(\mathop{\prod\limits_{l\in\mathrm{J}_0\backslash\{i,{i'}\}}} \limits_{\alpha_{l'}-\alpha_{l}\equiv0~(\mathrm{mod}~p)}x^{(\alpha_i\varepsilon_i)} x^{(\alpha_{^{l}}\varepsilon_l)})x^{((\overline{\alpha_i-1})\varepsilon_{i'})}x^{\bar{u}}).
\end{align*}
By computing the equation, we find that $a(\alpha,\bar{u})=0$ if $\alpha_{l}>0$ for $l\in\mathrm{J}_0\backslash\{i,i^{\prime}\}$ or $|\bar{u}|>0$.
Then we can suppose that
$$\phi(x_{l}x_{l'},x_{i})=\sum\limits_{0\leq\alpha\leq\pi}a(\alpha)x^{((\overline{\alpha_i-1})\varepsilon_{i'})}x^{(\alpha_i\varepsilon_i)}.$$
%Putting $k\in \mathrm{I}_1$
%\begin{align*}
%0=&\langle x_{k},\phi(x_{i}x_{i'},x_{i}) \rangle\\
%=&(D_{K}(x_{k})( \phi(x_{i}x_{i'},x_{i}))-G_{m}(x_{k})( \phi(x_{i}x_{i'},x_{i}))\\
%=&D_{k^{'}}( \phi(x_{i}x_{i'},x_{i}))\\
%=&D_{k^{'}}(\sum\limits_{0\leq\alpha\leq\pi} a(\alpha,i)(\mathop{\prod\limits_{l\in\mathrm{I}_0\backslash\{i,i^{'},m\}}} \limits_{\alpha_{l^{'}}-\alpha_{l}\equiv0(modp)} x^{(\alpha_{^{l}}\varepsilon_l)})x^{((\overline{\alpha_i-1})\varepsilon_{i'})}x^{(\alpha_i\varepsilon_i)}x^{\bar{u}} \rangle)
%\end{align*}
%If $\bar{u}$ has $k$ and $k^{'}$ in it, we have that $a(\alpha,i)=0$ if $|\bar{u}|>0$
Since $\mathrm{p}(x_{i}x_{i'})+\mathrm{p}(x_{i})=\overline{0}$ for any $ i\in \mathrm{I}_0$, by Lemma \ref{б³¬Ë«2}, we have
\begin{align*}
0=&\langle \phi(x_{i}x_{i'},x_{i}),\langle x_{i}x_{i'},x_{i}\rangle \rangle\\
=&\langle \phi(x_{i}x_{i'},x_{i}),-x_{i}\rangle\\
=&D_{K}(x_{i})( \phi(x_{i}x_{i'},x_{i}))-2\partial_{m}(x_{i})( \phi(x_{i}x_{i'},x_{i}))\\
=&\partial_{i'}( \phi(x_{i}x_{i'},x_{i}))\\
=&\partial_{i'}(\sum\limits_{0\leq\alpha\leq\pi} a(\alpha)x^{((\overline{\alpha_i-1})\varepsilon_{i'})}x^{(\alpha_i\varepsilon_i)}).
\end{align*}\\
By computing the equation, we find that $a(\alpha)=0$ if $\overline{\alpha_i-1}>0$.
Hence we get that $\alpha_i=1$.
Let $\lambda_i=a(\varepsilon_i)$.
From what has been discussed above, for any $i\in \mathrm{J}_0$ we have that
$$\phi(x_{l}x_{l'},x_{i})=-x_{i}=\lambda_{i}\langle x_{l}x_{l'},x_{i} \rangle,$$
where $\lambda_i$ is dependent on the second component.

%By Lemma \ref{WÉú³ÉÔªµÄȨ¿Õ¼ä}~(3), we can suppose that
%\begin{align*}
%\phi(x_{i}x_{i'},x_{i})=&\sum\limits_{0\leq\alpha\leq\pi}a(\alpha,i)(\mathop{\prod\limits_{i\in\mathrm{I}_0}} \limits_{\alpha_{i^{'}}-\alpha_{i}\equiv0(modp)} x^{(\alpha_{i}\varepsilon_i)})x^{\bar{u}};
%\end{align*}
%where $a(\alpha,i)\in \mathbb{F}$. For $k\in \mathrm{I}\backslash\{i,m\}$, we have that
%\begin{align*}
%0=&(-1)^{(\mathrm{p}(\phi)+\mathrm{p}(x_{i}x_{i'}))\mathrm{p}(x_{k})}(\phi(x_{i}x_{i'}, \langle x_{k},x_{i} \rangle)-\langle \phi(x_{i}x_{i'},x_{k}),x_{i} \rangle)\\
%=&\langle x_{k},\phi(x_{i}x_{i'},x_{i}) \rangle\\
%=&\langle x_{k},\sum\limits_{0\leq\alpha\leq\pi} a(\alpha,i)(\mathop{\prod\limits_{i\in\mathrm{I}_0\backslash\{m\}}} \limits_{\alpha_{i^{'}}-\alpha_{i}\equiv0(modp)} x^{(\alpha_{i}\varepsilon_i)})x^{(\alpha_m\varepsilon_m)}x^{\bar{u}} \rangle.
%\end{align*}
%When $k\in \mathrm{I}_1$. If $\bar{u}$ has k in it, we can deduce $a(\alpha,\bar{u})=0$;\\
%When $k\in \mathrm{I}_0 \backslash\{m\}$. $\alpha_{i}=1.$\\
%And we've seen that before $\alpha^{m}=0$.\\
%To sum up, for any $i\in \mathrm{I}_0$, we can get \\
%$$\phi(x_{l}x_{l'},x_{i})=a(0,i)(x_{i})=-a(0,i)(-x_{i})=\lambda_{i}\langle x_{i}x_{i'},x_{i} \rangle$$\\
%$$\phi(x_{i}x_{i'},x^{(q_i\varepsilon_i)})=\lambda_i\langle x_{i}x_{i'},x^{(q_i\varepsilon_i)} \rangle.$$
Similarly, we choose a fixed element $j\in \mathrm{I}_1$.
By Lemma \ref{WÉú³ÉÔªµÄȨ¿Õ¼ä} (3), we can suppose that
\begin{align*}
\phi(x_{j}x_{j'},x_{j})=&\sum\limits_{0\leq\alpha\leq\pi,~\bar{u}\in\mathbb{B}}a(\alpha,\bar{u},j)(\mathop{\prod\limits_{l\in\mathrm{J}_0}} \limits_{\alpha_{l'}-\alpha_{l}\equiv0~(\mathrm{mod}~p)} x^{(\alpha_{l}\varepsilon_l)})x^{(\alpha_m\varepsilon_m)}x^{\bar{u}}x_{j}
\end{align*}
where $a(\alpha,\bar{u},j)\in \mathbb{F}$.
By the definition of the shew-symmetric super-biderivation, we have
\begin{align*}
0=&\phi(x_{j}x_{j'},\langle 1,x_{j} \rangle)-\langle \phi(x_{j}x_{j'},1),x_{j}\rangle\\
=&(-1)^{(\mathrm{p}(\phi)+\mathrm{p}(x_{j}x_{j'}))\mathrm{p}(1)}\langle 1,\phi(x_{j}x_{j'},x_{j}) \rangle\\
%=&\langle 1,\phi(x_{i}x_{i'},x_{i}) \rangle\\
=&D_{K}(1)( \phi(x_{j}x_{j'},x_{j}))-2\partial_{m}(1)( \phi(x_{j}x_{j'},x_{j}))\\
=&2\partial_{m}( \phi(x_{j}x_{j'},x_{j}))\\
=&2\partial_{m}(\sum\limits_{0\leq\alpha\leq\pi,~\bar{u}\in\mathbb{B}}a(\alpha,\bar{u},j)
(\mathop{\prod\limits_{l\in\mathrm{J}_0}} \limits_{\alpha_{l'}-\alpha_{l}\equiv0~(\mathrm{mod}~p)} x^{(\alpha_{l}\varepsilon_l)})x^{(\alpha_m\varepsilon_m)}x^{\bar{u}}x_{j})\\
=&\sum\limits_{0\leq\alpha\leq\pi,~\bar{u}\in\mathbb{B}}a(\alpha,\bar{u},j)
(\mathop{\prod\limits_{l\in\mathrm{J}_0}} \limits_{\alpha_{l'}-\alpha_{l}\equiv0~(\mathrm{mod}~p)} x^{(\alpha_{l}\varepsilon_l)})x^{((\alpha_m-1)\varepsilon_m)}x^{\bar{u}}x_{j}.
\end{align*}
By computing the equation, we find that $a(\alpha,\bar{u},j)=0$ if $\alpha_{m}>0$.
For $k\in \mathrm{J}\backslash\{j,j^{\prime}\}$,
it is obvious that
\begin{align*}
0=&(-1)^{(\mathrm{p}(\phi)+\mathrm{p}(x_{l}x_{l'}))\mathrm{p}(x_{k})}(\phi(x_{l}x_{l'},\langle x_{k},x_{j} \rangle)-\langle \phi(x_{l}x_{l'},x_{k}),x_{j} \rangle)\\
=&\langle x_{k},\phi(x_{l}x_{l'},x_{j}) \rangle\\
=&\langle x_{k},\sum\limits_{0\leq\alpha\leq\pi,~\overline{u}\in\mathbb{B}} a(\alpha,\bar{u},j)(\mathop{\prod\limits_{l\in\mathrm{J}_0}} \limits_{\alpha_{l'}-\alpha_{l}\equiv0~(\mathrm{mod}~p)} x^{(\alpha_{l}\varepsilon_l)})x^{\bar{u}}x_{j} \rangle.
\end{align*}
Putting $k\in \mathrm{I}_1$, we can deduce $a(\alpha,\bar{u},j)=0$ if $|\bar{u}|>0$.
Putting $k\in\mathrm{J}_0$, we have that $a(\alpha,\bar{u},j)=0$ if $\alpha_{k}>0$.
Let $a(0,0,j)=\lambda_j$.
Hence, for any $j\in \mathrm{I}_1$, we have that $$\phi(x_{l}x_{l'},x_{j})=-a(0,0,j)x_{j}=\lambda_{j}\langle x_{l}x_{l'},x_{j} \rangle.$$
%
%We conclude that for any $i\in \mathrm{I}\setminus\{m\}$ and for any $x_{l}x_{l'}\in T$
%$$\phi(x_{l}x_{l'},x_{i})=-x_{i}=\lambda_{i}\langle x_{l}x_{l'},x_{i} \rangle,$$
%Set $\lambda_i:=-a(0,i)$ for $i\in\mathrm{I}$.\\
%where $\lambda_i$ is dependent on the second component.
The proof is completed.
\epf

\blem\label{3.11}
Suppose that $\phi$ is a  $\mathbb{Z}_{2}$-homogenous skew-symmetric super-biderivation on $K$. For any $x^{(q_i\varepsilon_i)}\in M$, where $1\leq q_i\leq \pi_i$, $i\in\mathrm{J}_0 $, there is an element $\lambda_{i}\in \mathbb{F}$ such that
$$\phi(x_{l}x_{l'},x^{(q_i\varepsilon_i)})=\lambda_{i}\langle x_{l}x_{l'},x^{(q_i\varepsilon_i)} \rangle.$$
\elem
\bpf
Without loss of generality,~we choose a fixed element $i\in \mathrm{I}_0$.
By Lemma \ref{Kб³¬Ë«1}, it is obvious that $\phi(x_{l}x_{l'},x^{(q_i\varepsilon_i)})=0$ for $l\in \mathrm{J}\setminus\{i,i'\}$ from $\langle x_{l}x_{l'},x^{(q_i\varepsilon_i)} \rangle=0$.
So we only need to consider the condition that $l=i$, it is clear that $\langle x_{i}x_{i'},x^{(q_i\varepsilon_i)} \rangle=-q_ix^{(q_i\varepsilon_i)}$.
If $q_i>1$, by Lemma \ref{б³¬Ë«1} and \ref{3.10}, we have
\begin{equation}
\begin{split}\label{²»ÖªµÀÊÇɶ1}
0=&\langle \phi(x_{k}x_{k'},x_{k}),\langle x_{i}x_{i'},x^{(q_i\varepsilon_i)} \rangle \rangle-(-1)^{\mathrm{p}(\phi)(\mathrm{p}(x_{k}x_{k'})+\mathrm{p}(x_{k}))}\langle \langle x_{k}x_{k'},x_{k} \rangle,\phi(x_{i}x_{i'},x^{(q_i\varepsilon_i)}) \rangle\\
=&\langle \lambda_k\langle x_{k}x_{k'},x_{k} \rangle,(-q_{i})x^{(q_i\varepsilon_i)} \rangle-\langle -x_{k},\phi(x_{i}x_{i'},x^{(q_i\varepsilon_i)}) \rangle\\
=&\langle x_{k},-\lambda_kq_{i}x^{(q_i\varepsilon_i)}+\phi(x_{i}x_{i'},x^{(q_i\varepsilon_i)}) \rangle.
\end{split}
\end{equation}
Because of $\mathcal{C}_{K_{-1}}(K)=\{f\in K|\langle f,x_{i} \rangle=0, \forall ~i\in\mathrm{I} \}= K_{-2}=\mathbb{F}1$, the equation (\ref{²»ÖªµÀÊÇɶ1}) implies that
$$\phi(x_{i}x_{i'},x^{(q_i\varepsilon_i)})=\lambda_i\langle x_{i}x_{i'},x^{(q_i\varepsilon_i)} \rangle+b,$$
where $\lambda_i$ is denoted in Lemma \ref{3.10} and $b\in\mathbb{F}$.
Since $\phi(x_{i}x_{i'},x^{(q_i\varepsilon_i)})\in K(m,n;\underline{t})_{(q_i\varepsilon_i)}$ by Lemma \ref{TµÄб˫³¬1}.
It is easily seen from Lemma \ref{WÉú³ÉÔªµÄȨ¿Õ¼ä} (2) that $K_{-2}\cap K_{(q_i\varepsilon_i)}=\emptyset$ for $q_i>1$.
So $b=0$ and $$\phi(x_{l}x_{l'},x^{(q_i\varepsilon_i)})=\lambda_i\langle x_{l}x_{l'},x^{(q_i\varepsilon_i)} \rangle,$$
where $\lambda_{i}$ is dependent on the second component.
\epf
\blem\label{3.12}
Suppose that $\phi$ is a $\mathbb{Z}_{2}$-homogenous skew-symmetric super-biderivation on $K$.
For any $x^{(q_m\varepsilon_m)}\in M$, where $0\leq q_{m}\leq \pi_{m}$, there is an element $\lambda\in \mathbb{F}$ such that
$$\phi(1,x^{(q_m\varepsilon_m)})=\lambda\langle 1,x^{(q_m\varepsilon_m)} \rangle.$$
\elem
\bpf
When $q_{m}=1$, we suppose that
$$\phi(1,x_m)=\sum_{0\leq \alpha\leq \pi,~u\in\mathbb{B}}c_{(\alpha,u)}x^{(\alpha)}x^{u},$$
where $c_{(\alpha,u)}\in \mathbb{F}$.
For $k\in\mathrm{J_{0}}$, by the definition of the skew-symmetric super-biderivation, we have the equation
\begin{align*}
0=&(-1)^{(p(\phi)+p(1))p(1)}(\phi(1,\langle1,x_m\rangle)-\langle \phi(1,1),x_m\rangle)\\
=&\langle 1,\phi(1,x_m)\rangle\\
=&\langle 1,\sum_{0\leq \alpha\leq \tau,~u\in\mathbb{B}}c_{\alpha}x^{(\alpha)}x^{u}\rangle\\
=&2\sum_{0<\alpha\leq \tau,~u\in\mathbb{B}}c_{\alpha}x^{(\alpha-\varepsilon_{m})}x^{u}.
\end{align*}
By computing the equation, we find that $c_{(\alpha,u)}=0$ if $\alpha_{m}-\varepsilon_{m}\geq0$.
We suppose that
$$\phi(1,x_m)=\sum_{0\leq \alpha_{\widehat{m}}\leq \tau_{\widehat{m}},~u\in\mathbb{B}}c_{\alpha_{(\widehat{m},u)}}x^{(\alpha_{\widehat{m}})}x^{u},$$
where $\widehat{m}$ represents an m-tuple with 0 as the m-th entry.
For $k\in\mathrm{J}$, by the definition of the skew-symmetric biderivation, we have the equation
\begin{align*}
0=&(-1)^{(p(\phi)+p(1))p(x_{m})}(\phi(x_{k},\langle1,x_m\rangle)-\langle \phi(x_k,1),x_m\rangle)\\
=&\langle x_{k},\phi(1,x_m)\rangle\\
=&\langle x_{k},\sum_{0\leq \alpha_{\widehat{m}}\leq \tau_{\widehat{m}},~u\in\mathbb{B}}c_{\alpha_{(\widehat{m},u)}}x^{(\alpha_{\widehat{m}})}x^{u}\rangle\\
=&\sum_{0\leq \alpha_{\widehat{m},~u\in\mathbb{B}}\leq \tau_{\widehat{m}}}c_{\alpha_{(\widehat{m},u)}}x^{(\alpha_{(\widehat{m},u)}-\varepsilon_{k'})}x^{u}.
\end{align*}
By computing the equation, we find that $c_{\alpha_{(\widehat{m},u)}}=0$ if $\alpha_{\widehat{m}}-\varepsilon_{k'}>0$ or $|u|>0$.
Then we can suppose that
$$\phi(1,x_{m})=c_{0}.$$
Set $\lambda=\frac{c_0}{2}$, then we can get that
$$\phi(1,x_m)=\lambda\langle1,x_m\rangle.$$

When $q_{m}\geq 2$, we suppose that $$\phi(1,x^{(q_{m}\varepsilon_m)})=\sum_{0\leq \alpha\leq \tau,~u\in\mathbb{B}}c_{(\alpha,u)}x^{(\alpha)}x^{u}.$$
By Lemma \ref{б³¬Ë«1} and the conclusion of the case $q_{m}=1$, we have the equation
\begin{align*}
0=&\langle \phi(1,x_{m}),\langle 1,x^{(q_{m}\varepsilon_{m})}\rangle\rangle-\langle\langle 1,x_{m}\rangle,\phi(1,x^{(q_{m}\varepsilon_{m})})\rangle\\
=&\langle\lambda\langle1,x_{m}\rangle,\langle 1,x^{(q_{m}\varepsilon_{m})}\rangle\rangle-\langle\langle 1,x_{m}\rangle,\sum_{0\leq \alpha\leq \tau,~u\in\mathbb{B}}c_{\alpha}x^{(\alpha)}x^{u}\rangle\\
=&4\lambda\langle1,x^{((q_{m-1})\varepsilon_{m})}\rangle-2\langle1,\sum_{0\leq \alpha\leq \tau,~u\in\mathbb{B}}c_{(\alpha,u)}x^{(\alpha)}x^{u}\rangle\\
=&8\lambda x^{((q_{m}-2)\varepsilon_{m})}-4\sum_{0\leq \alpha\leq \tau,~u\in\mathbb{B}}c_{(\alpha,u)}x^{(\alpha-\varepsilon_{m})}x^{u}.
\end{align*}
By computing the equation, we find that $c_{(\alpha,u)}=0$ if $\alpha-(q_{m}-1)\varepsilon_{m}\neq0$ or $|u|>0$. And $c_{(q_{m}-1)}=2\lambda$.
We suppose that
$$\phi(1,x^{(q_{m}\varepsilon_m)})=\sum_{0\leq \alpha_{\widehat{m}}\leq \tau_{\widehat{m}}}c_{\alpha_{\widehat{m}}}x^{(\alpha_{\widehat{m}})}+2\lambda x^{((q_{m}-1)\varepsilon_m)}.$$
For any $i \in \mathrm{J_{0}}$, by Lemma \ref{б³¬Ë«1}, we have the equation
\begin{equation*}
\begin{split}
0=&\langle\phi(1,x^{(q_{m}\varepsilon_{m})}),\langle x_m,x_i\rangle\rangle-\langle\langle1,x^{q_{m\varepsilon_{m}}}\rangle,\phi(x_m,x_i)\rangle\\
=&\langle\sum_{0\leq \alpha_{\widehat{m}}\leq \tau_{\widehat{m}}}c_{\alpha_{\widehat{m}}}x^{(\alpha_{\widehat{m}})}+2\lambda x^{((q_{m}-1)\varepsilon_m)},-x_i\rangle
-\langle2x^{((q_{m}-1)\varepsilon_n)},\lambda_i\langle x_m,x_i \rangle\rangle\\
=&\langle x_i,\sum_{0\leq \alpha_{\widehat{m}}\leq \tau_{\widehat{m}}}c_{\alpha_{\widehat{m}}}x^{(\alpha_{\widehat{m}})}+2\lambda x^{((q_{m}-1)\varepsilon_m)}-2\lambda_ix^{((q_{m}-1)\varepsilon_m)} \rangle.
\end{split}
\end{equation*}
Since $\mathcal{C}_{K_{-1}}(K)= K_{-2}$, we have that
$$\sum_{0\leq \alpha_{\widehat{m}}\leq \tau_{\widehat{m}}}c_{\alpha_{\widehat{m}}}x^{(\alpha_{\widehat{m}})}+2\lambda x^{((q_{m}-1)\varepsilon_m)}-2\lambda_ix^{((q_{m}-1)\varepsilon_m)}\in\mathbb{F}.$$
Then we have $c_{\alpha_{\widehat{m}}}=0$ for $\alpha_{\widehat{m}}>0$ and $\lambda=\lambda_i$ for $i \in \mathrm{J_{0}}$.
Then we can get that $$\phi(1,x^{(q_{m}\varepsilon_m)})=c_0+2\lambda x^{((q_{m}-1)\varepsilon_m)}.$$
Utilizing the definition of the skew-symmetry biderivation, by Lemma \ref{Kб³¬Ë«1}, we have that
\begin{align*}
0=&\langle\phi(1,x^{q_{m}\varepsilon_{m}}),\langle1,x^{q_{m}\varepsilon_{m}}\rangle\rangle\\
=&\langle c_0+2\lambda x^{((q_{m}-1)\varepsilon_m)},x^{(q_{m}-1)\varepsilon_{m}}\rangle\\
=&2c_0x^{(q_{m}-2)\varepsilon_{m}}.
\end{align*}
It is obvious that $c_0=0$ for $p>2$.
So we can get that
$$\phi(1,x^{(q_{m}\varepsilon_m)})=\lambda\langle1,x^{(q_{m}\varepsilon_m)}\rangle.$$
The proof is complete.
\epf

\begin{re}\label{3.13}
We claim that $\lambda_1=\cdots=\lambda_m=\cdots=\lambda_{m+n}$.
Choose two mutually different elements $i,j\in \mathrm{J}$.
Since the characteristic $p>3$, there are two positive integers $q_i$ and $q_m$, which are greater than 1 and are neither congruent to 0 modulo $p$, such that we have
\begin{align*}
0=&\langle \phi(x_{i}x_{i'},x_{i}),\langle 1,x^{(q_m\varepsilon_m)}\rangle\rangle-\langle\langle x_{i}x_{i'},x_{i}\rangle,\phi(1,x^{(q_m\varepsilon_m)})\rangle\\
=&\langle\lambda_{i}\langle x_{i}x_{i'},x_{i}\rangle,2x^{((q_m-1)\varepsilon_m)}\rangle-\langle -x_{i},\lambda\langle 1,x^{(q_m\varepsilon_m)}\rangle\rangle\\
=&\langle -\lambda_{i}x_{i},2x^{((q_m-1)\varepsilon_m)}\rangle-\langle -x_{i},2\lambda x^{((q_m-1)\varepsilon_m)}\rangle\\
=&-2\lambda_{i}\langle x_{i},x^{((q_m-1)\varepsilon_m)}\rangle+2\lambda\langle x_{i},x^{((q_m-1)\varepsilon_m)}\rangle\\
=&2(\lambda-\lambda_{i})\langle x_{i},x^{((q_m-1)\varepsilon_m)}\rangle\\
=&2(\lambda-\lambda_{i})x^{(\varepsilon_{i}+(q_m-1)\varepsilon_m)}.
\end{align*}
By direct calculation,
it is easily seen that $\lambda_i=\lambda$ for any $i\in\mathrm{ I }$.
Set $\lambda:=\lambda_1=\cdots=\lambda_m=\lambda_{m+1}=\cdots=\lambda_{m+n}$.
Then we can conclude that for any $x^{(q_i\varepsilon_{i})}\in M$, $1\leq q_i\leq \pi_i$ and $x_{l}x_{l'}\in T$, there is an element $\lambda\in \mathbb{F}$ such that
$$\phi(x_{l}x_{l'},x^{(q_i\varepsilon_{i})})=\lambda\langle x_{l}x_{l'},x^{(q_i\varepsilon_i)}\rangle,$$
where $\lambda$ depends on neither $x^{(q_i\varepsilon_i)}$ nor $x_{l}x_{l'}$.
\end{re}

\bthm
Let $K$ be the contact Lie superalgebra $K(m,n;\underline{t})$ over the prime field $\mathbb{F}$ of the characteristic $p>3$, where $m,~n\in \mathbb{N}+1$ and $\underline{t}=(t_1,t_2,\ldots,t_m)$ is an $m$-tuple of positive integers.
Then
$$\mathrm{BDer}(K)=\mathrm{IBDer}(K).$$
\ethm
\bpf
Suppose that $\phi$ is a skew-symmetric super-biderivation on $K$.
By Lemmas \ref{3.09} and \ref{3.10}, there is an element $\lambda\in\mathbb{F}$ such that
$\phi(x_{i}x_{i'},x_{i})=\lambda\langle x_{i}x_{i'},x_{i}\rangle$ for all $i\in \mathrm{J}$.
For any $x^{(\alpha)}x^u, x^{(\beta)}x^v\in K$ and $x_{l}x_{l'}\in T_{k}$, by Lemma \ref{б³¬Ë«1} and Remark \ref{б³¬Ë«ÊÇżµÄ}, we have the equation
\begin{equation*}
\begin{split}
0&=\langle \phi(x_{l}x_{l'},x_{l}),\langle x^{(\alpha)}x^u,x^{(\beta)}x^v \rangle \rangle-\langle \langle x_{l}x_{l'},x_{l} \rangle,\phi(x^{(\alpha)}x^u,x^{(\beta)}x^v) \rangle\\
&=\langle \langle x_{l}x_{l'},x_{l} \rangle,\lambda\langle x^{(\alpha)}x^u,x^{(\beta)}x^v \rangle \rangle-\langle \langle x_{l}x_{l'},x_{l} \rangle,\phi(x^{(\alpha)}x^u,x^{(\beta)}x^v) \rangle\\
&=\langle x_{l},\phi(x^{(\alpha)}x^u,x^{(\beta)}x^v)-\lambda\langle x^{(\alpha)}x^u,x^{(\beta)}x^v \rangle \rangle.
\end{split}
\end{equation*}
Since $\mathcal{C}_{K_{-1}}(K)= K_{-2}$, we have that
$$\phi(x^{(\alpha)}x^u,x^{(\beta)}x^v)=\lambda\langle x^{(\alpha)}x^u,x^{(\beta)}x^v \rangle+b,$$
where $\lambda$ is denoted in Remark \ref{3.13} and $b\in\mathbb{F}$.
By Lemma \ref{б³¬Ë«1} and Remark \ref{б³¬Ë«ÊÇżµÄ}, we have
\begin{equation*}
\begin{split}
0&=\langle\phi(x^{(\alpha)}x^u,x^{(\beta)}x^v),\langle 1,x^{(2\varepsilon_{m})}\rangle\rangle-\langle \langle x^{(\alpha)}x^u,x^{(\beta)}x^v \rangle,\phi(1,x^{(2\varepsilon_{m})}) \rangle\\
&=\langle \lambda\langle x^{(\alpha)}x^u,x^{(\beta)}x^v \rangle+b,2x_{m}\rangle-\langle \langle x^{(\alpha)}x^u,x^{(\beta)}x^v \rangle,\lambda\langle 1,x^{(2\varepsilon_{m})} \rangle \rangle\\
&=\langle \lambda\langle x^{(\alpha)}x^u,x^{(\beta)}x^v \rangle+b,2x_{m}\rangle-\langle \langle x^{(\alpha)}x^u,x^{(\beta)}x^v \rangle,\lambda2x_{m}\rangle\\
&=\langle b,2x_{m}\rangle\\
&=4b.
\end{split}
\end{equation*}
Then $b=0$. %for all $s\in\mathrm{I}$.
Hence, $\phi(x^{(\alpha)}x^u,x^{(\beta)}x^v)=\lambda\langle x^{(\alpha)}x^u,x^{(\beta)}x^v \rangle$ ~for any $x^{(\alpha)}x^u,x^{(\beta)}x^v\in K $ and $\phi$ is an inner super-biderivation.
\epf

\noindent {\bf Acknowledgements}\quad The authors would like to
thank the referee for valuable comments and suggestions on this
article.


\begin{thebibliography}{99}
%1
\bibitem{Bm} M. Bre$\rm{\check{s}}$ar, Commuting maps: a survey, Taiwanese J. Math, 8 (2004), 361--397.
%2
\bibitem{WdYxCHz} D. Wang, X. Yu, Z. Chen, Biderivations of the parabolic subalgebras of simple Lie algebras, Comm. Algebra, 39 (2011), 4097--4104.
%3
\bibitem{CHz} Z. Chen, Biderivations and linear commuting maps on simple generalized Witt algebras over a field, Electron. J. Linear Algebra, 31 (2016), 1--12.
%4
\bibitem{HxWdXc} X. Han, D. Wang, C. Xia, Linear commuting maps and biderivations on the Lie algebras $\mathcal{W}$(a,b), J. Lie Theory, 26 (2016), 777--786.
%5
\bibitem{Tx}X. Tang, Biderivations of finite-dimensional complex simple Lie algebras, Linear Multilinear Algebra, 66  (2018), 250--259.
%6
\bibitem{WdYx} D. Wang, X. Yu, Biderivations and linear commuting maps on the  Schr$\rm{\ddot{o}}$dinger-Virasoro Lie algebra, Comm. Algebra, 41 (2013), 2166--2173.
%7
\bibitem{BmZHk} M. Bre$\rm{\check{s}}$ar, K. Zhao, Biderivations and commuting linear maps on Lie algebras, J. Lie Theory, 28 (2018), 885--900.
%8
\bibitem{CHyCHly} Y. Chang, L. Chen, Biderivations and linear commuting maps on the restricted Cartan-type Lie algebras $W(n;\underline{1})$ and $S(n;\underline{1})$, Linear Multilinear Algebra, DOI:10.1080/03081087.2018.1465525.
%9
\bibitem{CHyCHlyZHx} Y. Chang, L. Chen, X. Zhou,  Biderivations and linear commuting maps on the restricted Cartan-type Lie algebras $H(n;\underline{1})$,
 Comm. Algebra 47 (2019), 1311-1326.
%10
\bibitem{FgDx} G. Fan, X. Dai, Super-biderivations of Lie superalgebras, Linear Multilinear Algebra, 65 (2017), 58--66.
%11
\bibitem{XcWdHx} C. Xia, D. Wang, X. Han, Linear super-commuting maps and super-biderivations on the super-Virasoro algebras, Comm. Algebra, 44 (2016), 5342--5350.
%12
\bibitem{YjTx} J. Yuan, X. Tang, Super-biderivations of classical simple Lie superalgebras, Aequationes Math., 92 (2018), 91--109.
%13
\bibitem{YchLchYc}Y. Chang, L. Chen, Y. Cao, Super-biderivations of the generalized Witt Lie superalgebra W(m,n;\underline{t}), Linear Multilinear Algebra, DOI:10.1080/03081087.2019.1593312.
%14
\bibitem{FmQzh}F. Ma, Q. Zhang, Derivation algebras for $K$-type modular Lie superalgebras,  J. Math. (Wuhan), 20 (2000), 431--435.
%15
\bibitem{BgWl}B. Guan, W. Liu, Derivations of the even part into the odd part for modular contact superalgebra,
J. Math. (Wuhan) 32 (2012), 402--414.
%16
\bibitem{BgLch}B. Guan, L. Chen, Derivations of the even part of contact Lie superalgebra, J. Pure Appl. Algebra 216 (2012), 1454--1466.
%17
\bibitem{YZhWl}Y. Zhang, W. Liu, Moduliar Lie superalgebras, Science Press, Beijing, 2005.
%18
\bibitem{YwYZh}Y. Wang, Y. Zhang, The associative forms of the graded Cartan type Lie superalgebras.  Adv. Math. (in Chinese), 29 (2000), 65--70.
%
%
%
%\bibitem{Bai&Liu} W. Bai and W. Liu, Superderivations for modular graded Lie superalgebras of Cartan-type, Algebr. Represent. Theory, 17 (2014), 69--86.
%%\bibitem{Bdb} D. Benkovi$\rm{\check{s}}$c, Biderivations of triangular algebras. Linear Algebra Appl. 431 (2009), 1587--1602.
%%\bibitem{Bdg} D. Benkovi$\rm{\check{s}}$c, Generalized Lie derivations on triangular algebras. Linear Algebra Appl, 434 (2011), 1532--1544.
%\bibitem{HxWL}H. Xu, L. Wang, The properties of biderivations on Heisenberg superalgebras. Math.
%    Aeterna 5(2) (2015), 285--291
%\bibitem{Ma&Zhang} F. Ma and Q. Zhang, Derivation algebra of modular Lie superalgebra $K$ of Cartan-type, J. Math. (Wuhan), 20(4)(2000), 431--435.
%\bibitem{NG} N. Ghosseiri, On biderivations of upper triangular matrix rings. Linear Algebra Appl.
%    438(1) (2013), 250--260.
%\bibitem{Wang&Zhang} Y. Wang and Y. Zhang, Derivation algebra $\mathbf{Der}(H)$ and central extensions of Lie superalgebras, Comm. Algebra, 32 (2004), 4117--4131.
%\bibitem{Zhy} Y. Zhang, Finite-dimensional Lie superalgebras of Cartan-type over fields of prime characteristic, Chinese Sci. Bull., 42 (1997), 720--724.
%\bibitem{Zhang&Zhang} Q. Zhang and Y. Zhang, Derivation algebras of the modular Lie superalgebras $W$ and $S$ of Cartan-type,  Acta Math. Sci. Ser. B (Engl. Ed.), 20 (2000), 137--144.
\end{thebibliography}
\end{document}